\documentclass{article}

\usepackage{arxiv}

\usepackage[utf8]{inputenc} 
\usepackage[T1]{fontenc}    
\usepackage{hyperref}       
\usepackage{url}            
\usepackage{booktabs}       
\usepackage{amsfonts}       
\usepackage{nicefrac}       
\usepackage{microtype}      

\usepackage{amsmath,amssymb,amsthm,listings,xcolor}
\newtheorem{assumption}{Assumption}
 
\colorlet{punct}{red!60!black}
\definecolor{background}{HTML}{EEEEEE}
\definecolor{delim}{RGB}{20,105,176}
\colorlet{numb}{magenta!60!black}
\lstdefinelanguage{json}{
    basicstyle=\normalfont\ttfamily,
    numbers=left,
    numberstyle=\scriptsize,
    stepnumber=1,
    numbersep=8pt,
    showstringspaces=false,
    breaklines=true,
    frame=lines,
    backgroundcolor=\color{background},
    literate=
     *{0}{{{\color{numb}0}}}{1}
      {1}{{{\color{numb}1}}}{1}
      {2}{{{\color{numb}2}}}{1}
      {3}{{{\color{numb}3}}}{1}
      {4}{{{\color{numb}4}}}{1}
      {5}{{{\color{numb}5}}}{1}
      {6}{{{\color{numb}6}}}{1}
      {7}{{{\color{numb}7}}}{1}
      {8}{{{\color{numb}8}}}{1}
      {9}{{{\color{numb}9}}}{1}
      {:}{{{\color{punct}{:}}}}{1}
      {,}{{{\color{punct}{,}}}}{1}
      {\{}{{{\color{delim}{\{}}}}{1}
      {\}}{{{\color{delim}{\}}}}}{1}
      {[}{{{\color{delim}{[}}}}{1}
      {]}{{{\color{delim}{]}}}}{1},
}
\definecolor{capri}{rgb}{0.0, 0.75, 1.0}

\usepackage{graphicx,wrapfig,multirow,diagbox,caption,subcaption,enumitem}
\usepackage{verbatim,algpseudocode}
\usepackage[boxruled]{algorithm2e}

\usepackage[numbers,sort,compress]{natbib}
\usepackage{bm,mathrsfs}

\newcommand{\ud}{\mathrm{d}}
\newcommand{\mr}{\mathbb{R}}

\DeclareMathOperator*{\st}{s.t.}

\title{Computational Graph Representation of Equations System Constructors in Hierarchical Circuit Simulation}
\author{Zichao Long
\And
Lin Li
\And
Lei Han
\And
Xianglong Meng
\And
Chongjun Ding
\And
Ruiyan Li
\And
Wu Jiang
\And
Fuchen Ding
\And
Jiaqing Yue
\And
Zhichao Li
\And
Yisheng Hu
\And
Ding Li
\And
Heng Liao
\thanks{Corresponding author}
\AND
\normalfont{HiSilicon}
}
\begin{document}
\maketitle
\begin{abstract}
  Equations system constructors of hierarchical circuits play a central role in
  device modeling, nonlinear equations solving, and circuit design automation.
  However, existing constructors present limitations in applications to different
  extents. For example, the costs of developing and reusing device models ---
  especially coarse-grained equivalent models of circuit modules ---
  remain high while parameter sensitivity analysis is complex and inefficient.
  Inspired by differentiable programming and
  leveraging the ecosystem benefits of open-source software, we propose
  an equations system constructor using the computational graph representation,
  along with its JSON format netlist, to address these limitations.
  This representation allows for runtime dependencies between signals and
  subcircuit/device parameters. The proposed method streamlines the model
  development process and facilitates end-to-end computation of gradients of
  equations remainders with respect to parameters.
  This paper discusses in detail the overarching concept of hierarchical
  subcircuit/device decomposition and nested invocation by drawing parallels to
  functions in programming languages, and introduces
  rules for parameters passing and gradient propagation across
  hierarchical circuit modules.
  The presented numerical examples, including (1) an uncoupled CMOS model
  representation using "equivalent circuit decomposition+dynamic parameters"
  and (2) operational amplifier (OpAmp) auto device sizing, have demonstrated
  that the proposed method supports circuit simulation and design and
  particularly subcircuit modeling with improved efficiency, simplicity,
  and decoupling compared to existing techniques.
\end{abstract}

\newpage
\tableofcontents

\newpage
\section{Introduction}
\begin{figure}[htpb]
	\centering
	\includegraphics[width=1.0\textwidth]{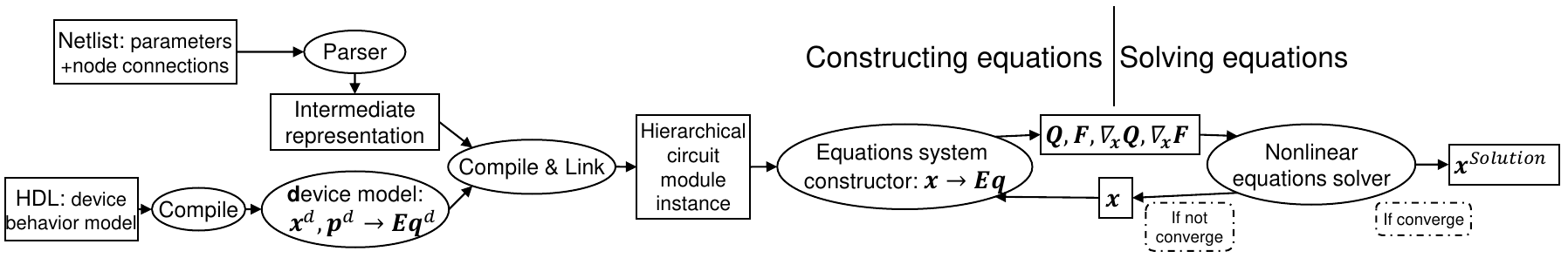}
	\caption{Simulator flowchart}
	\label{fig:simulator-flowchart}
\end{figure}

Analog circuit simulation (Figure \ref{fig:simulator-flowchart}) \cite{nagel1971computer,mccalla1971bias,Nagel:M382} has become one of the most important tools in the analog circuit EDA toolchain to assist verification and design. This is due to long-term research accumulation in algebraic differential equation theory \cite{estevez2000structural,kunkel2006differential,gunther2005modelling}, device modeling \cite{chauhan2012bsim,ezaki2008physics,gildenblat2006psp}, equations system construction \cite{hachtel1971sparse,ho1975modified,fijnvandraat2002time}, nonlinear equations solvers \cite{nastov2007fundamentals,najm2010circuit}, hardware description languages (HDLs) \cite{verilog2014verilog,ieee2006ieee-1364-2005,lemaitre2002adms,christen1999vhdl,pecheux2005vhdl}, and more.
However, we still see inconveniences with popular HDLs and simulators in reusing coarse-grained circuit module behavior models, introducing multi-physical effects, and applying gradient optimization methods to automatic design.
\begin{itemize}[partopsep=0pt,topsep=0pt,itemsep=0pt,parsep=0pt]
	\item
	For device and circuit modeling, HDL technology is often used to create a device behavior model, which is then compiled \cite{lemaitre2002adms} into a program that the simulator can invoke. In particular, multi-port circuit behavior modeling often goes through two steps: function fitting and HDL implementation, for example, neural network fitting \cite{meijer2001neural,zhang2017artificial}+HDL or Volterra polynomial \cite{zhang2014large}+SPICE netlists \cite{zhang2019new,roymohapatra2019novel}.
	\item
	To achieve an optimal combination of circuit parameters in the automation of analog circuit design \cite{razavi2002design,silveira1996g,jespers2009gm}, it is necessary to first optimize the device size. Prior to 2010, researchers developed solutions based on gradient optimization using parameter sensitivity analysis \cite{zhan2004optimization,agrawal2006circuit}. To tune design variables, modern process and software technologies require that software uses callbacks to modify model parameters. However, such an approach hinders users from obtaining gradient information of the variables. This has pushed recent research to shift its focus toward black-box methods, such as local sampling for gradient reconstruction \cite{huang2013efficient,nieuwoudt2005multi,peng2016efficient}, surrogate models \cite{girardi2011analog,lyu2018batch,wang2014enabling,lyu2017efficient}, and reinforcement learning \cite{tang2018parametric}.
\end{itemize}
Many works attempt to address the above inconveniences from both the model compilation and the gradient acquisition perspectives. For example, \citet{mahmutoglu2018new} developed the Verilog-AMS compiler that can run using Matlab/Octave, \citet{kuthe2020verilogae} can obtain more information about the equations and internal derivatives when the Verilog-AMS circuit module is compiled, and \citet{hu2020adjoint} provides an efficient implementation of adjoint equations for transient simulation. However, none of these works explore functional support possibilities from the perspective of systems equations construction.

Indeed, one significant cause of these inconveniences is that the analog HDLs contain many complex and even bloated features, which are necessary for these HDLs to simultaneously support the representation of structural and behavioral information. Such features include the automatic differentiation required for analog simulation, and polynomial interpolation that may be used in modeling \cite[Sections 4.5.6, 9.21]{verilog2014verilog}. The intertwining of structural and behavioral information has isolated analog HDLs from the open ecosystems of other programming languages and tools, and also created high barriers for developing analog EDA tools. Furthermore, only static circuit parameters independent of signal values can be passed between nested circuit modules and simulation runtime variables are allowed only within modules (\cite[Sections 3.4, 6]{verilog2014verilog},\cite[Section 4.10]{ieee2006ieee-1364-2005}) --- this is not conducive to reuse and development of coarse-grained circuit models.

\begin{figure}[htpb]
	\centering
	\includegraphics[width=0.7\textwidth]{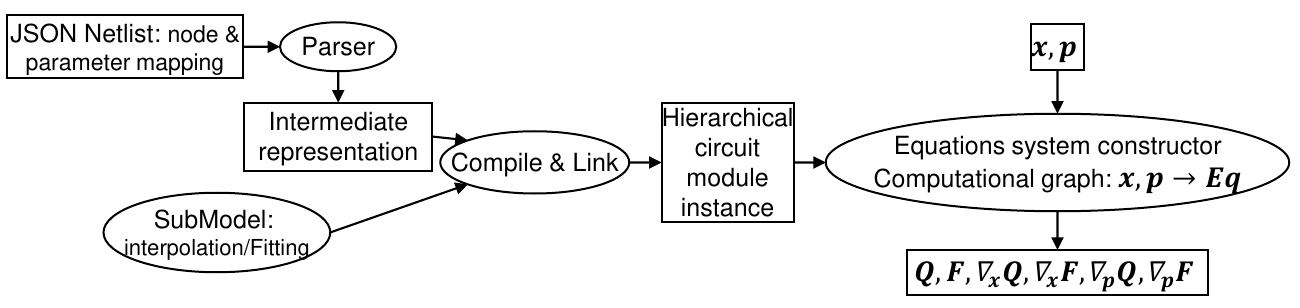}
	\caption{Flowchart: create computational graph from a JSON netlist}
	\label{fig:flowchart}
\end{figure}

Thriving technologies such as deep learning \cite{goodfellow2016deep} and automatic differential programming \cite{baydin2018automatic} have contributed to a number of research fields in scientific computing, including data-driven multi-scale modeling and inverse problems \cite{zhang2018deep,long2018pde,long2019pde}. Inspired by such works, this paper proposes a computational graph implementation of an equations system constructor that works in hierarchical circuit simulation \cite{fijnvandraat2002time,ter1999numerical,mukherjee1999hierarchical,tcherniaev2003transistor}, along with the corresponding JSON netlist and compilation method (Figure \ref{fig:flowchart}). The internal and external variables, design variables, model input parameters, and corresponding gradients of circuit modules are processed in a unified manner.
This work takes circuit modules as the basic units of a computational graph and supports decoupled representation of circuit models using "equivalent circuit decomposition using JSON netlist+submodel for computing dynamic parameters". The advantages of this approach are twofold: (1) JSON format is easy to parse, and its basic data types "dictionary+list" are sufficient to represent circuit structure information. And (2), a submodel can be implemented with the help of the automatic differentiation capability of Julia\cite{Bezanson_Julia_A_fresh_2017}. Based on these two advantages, the proposed method in this paper simplifies circuit modeling and enhances the gradient acquisition capability of simulation tools.

Section \ref{sec:Joanna} describes the processing of structural information related to circuit modules in a computational graph. Similar to defining, compiling, representing, and executing functions in a programming language \cite{aho2007compilers,muchnick1997advanced,appel2004modern}, the process involves the following steps:
module definition (netlist), parsing and compilation, data structures of module instances, and graph executor (Figure \ref{fig:equations-system-constructor}). Additionally, Section \ref{subsec:EvalCompositeSubCkt} introduces how to define and use SubModel for computing intrinsic dynamic parameters to provide behavioral information of each circuit module. Section \ref{sec:applications} presents two application examples, one for device modeling and the other for a joint solution of DC/AC analysis and device sizing under different process, voltage, and temperature (PVT) combinations.

\section[Hierarchical Circuit Equations System Constructor: Computational Graph]{Hierarchical Circuit Equations System Constructor: \\ Computational Graph}\label{sec:Joanna}
\begin{figure}[htpb]
	\centering
	\begin{subfigure}{0.49\textwidth}
		\includegraphics[width=\textwidth]{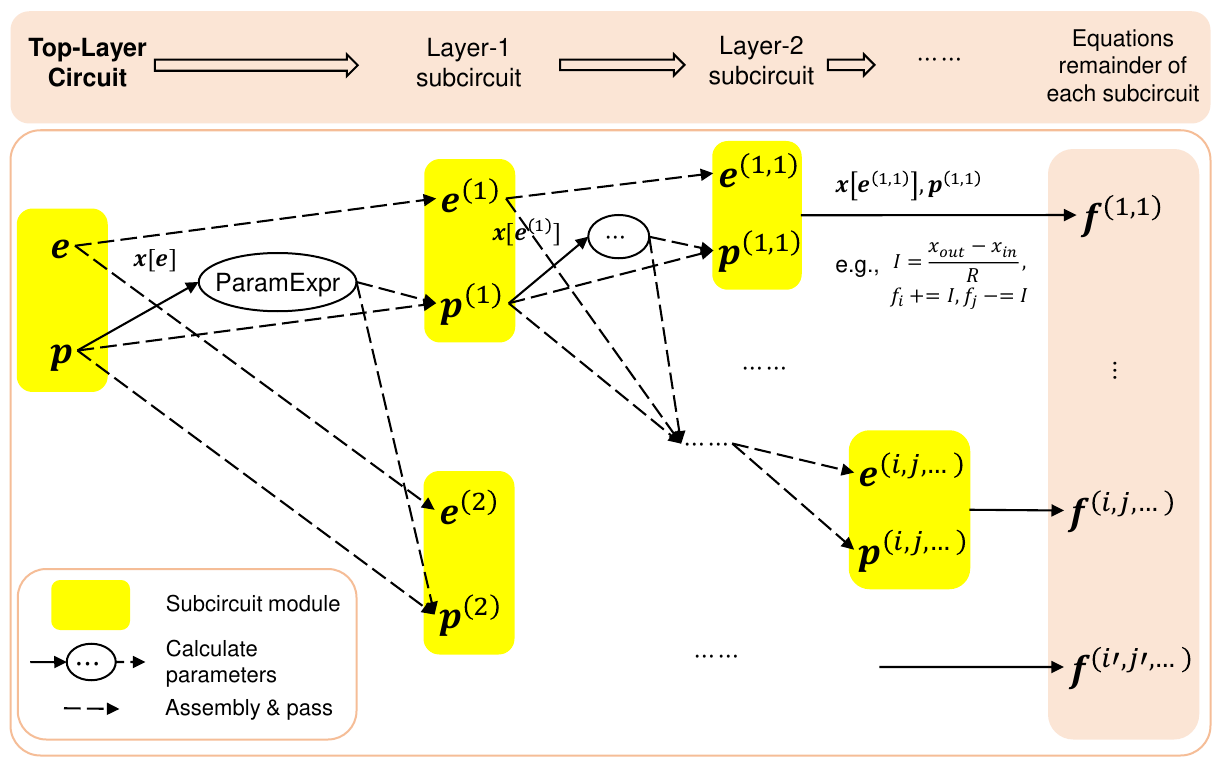}
		\caption{Existing method using static parameters}
		\label{fig:static-engine}
	\end{subfigure}
	\begin{subfigure}{0.49\textwidth}
		\includegraphics[width=\textwidth]{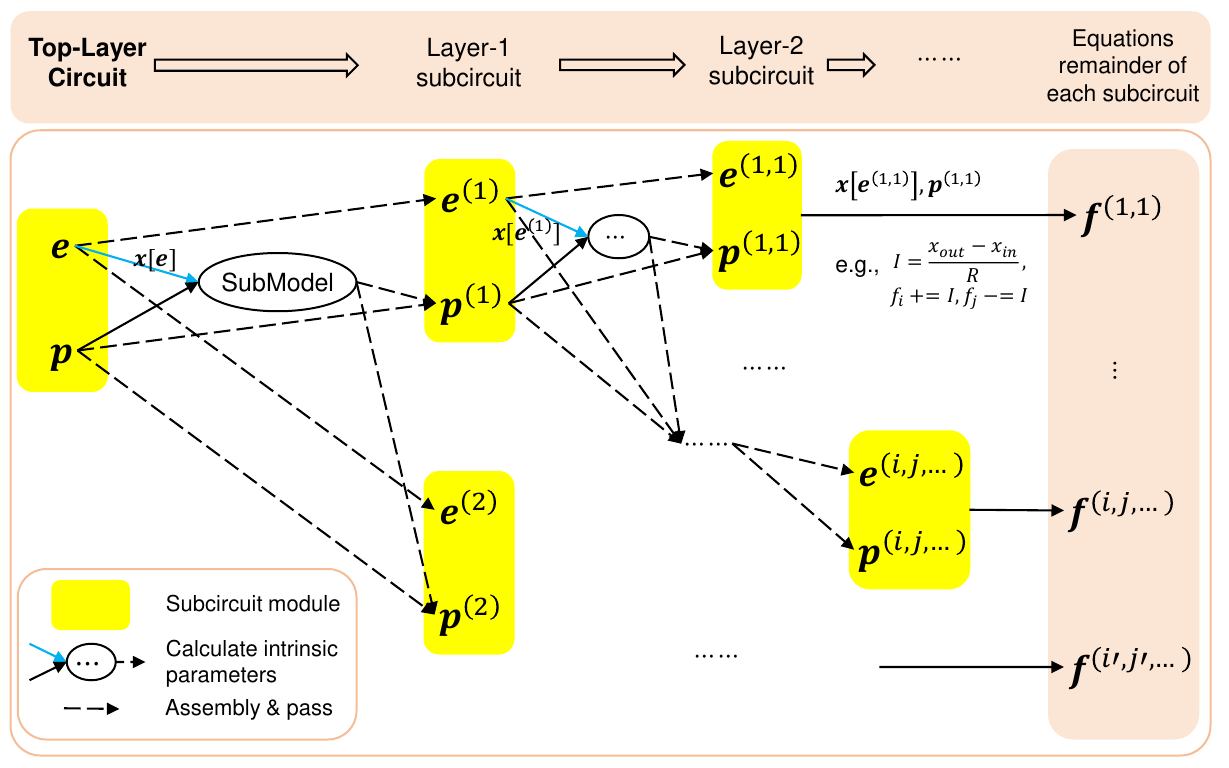}
		\caption{Computational graph using dynamic parameters}
		\label{fig:computational-graph}
	\end{subfigure}
	\caption{Equations system constructor for hierarchical circuits: existing method v.s. computational graph. Each computing unit corresponds to a subcircuit, which can be decomposed into smaller subcircuits --- the smallest granularity subcircuit are the "basic elements". $\bm{e}^{(\cdots)}$ denotes the internal and external nodes of a subcircuit, $\bm{x}$ a generalized signal, such as node bias or branch current, $\bm{p}^{(\cdots)}$ the input parameters of a subcircuit, for example, the device size (alternatively, a non-linear capacitance, inductance, or current value of a basic element, functioning as a dynamic parameter), and $\bm{f}^{(\cdots)}$ the contribution of a subcircuit to the remainder of the simulation equations. Calculation of the Jacobian matrix is omitted in this schematic. In the existing method, a lower-layer parameter of a circuit module does not depend on the signal value ($\bm{x}[\bm{e}]$), and calculation of ParamExpr can be completed when the netlist is built. In the computational graph method, the dynamic parameters are derived from the {\color{capri}intrinsic parameters} output by SubModel --- these parameters are obtained at the calculation runtime of the systems of equations.}
	\label{fig:equations-system-constructor}
\end{figure}
For any $N$ generalized signals (i.e., equation unknowns) $\bm{x}\in\mr^N$ and $M$ signal-independent input variables or parameters $\bm{p}\in\mr^M$, the mathematical meaning of circuit simulation is constructing and solving the following conservation algebraic differential equations \cite{najm2010circuit,gunther2005modelling,hu2020adjoint}
\begin{equation}\label{eq:flat-equation}
\bm{f}(\bm{\dot{x}}(t), \bm{x}(t), \bm{p})
\triangleq\frac{\ud\bm{Q}(\bm{x},\bm{p})}{\ud t}+\bm{F}(\bm{x},\bm{p})
=\bm{0},\tag{Eq (Flat)}
\end{equation}
where, $\bm{Q}$ denotes the dynamic part of the remainder of the equation (e.g., the charge of different capacitors and magnetic flux of different inductors) and $\bm{F}$ denotes the static part of the remainder of the equation (e.g., the total input DC current at each node and the voltage drop equation). Modern simulators tend to decompose and construct \ref{eq:flat-equation} based on circuit hierarchy, which are easier to understand and parallelize:
\begin{equation}\label{eq:hierarchical-equation}
\begin{split}
\bm{f}(\bm{\dot{x}},\bm{x},\bm{p})
& = \bm{f}^{(1)}(\bm{\dot{x}}^{(1)},\bm{x}^{(1)},\bm{p}^{(1)})
+ \bm{f}^{(2)}(\bm{\dot{x}}^{(2)},\bm{x}^{(2)},\bm{p}^{(2)})
+ \cdots \\
& = \bm{f}^{(1,1)}+\bm{f}^{(1,2)}+\cdots+\bm{f}^{(2,1)}+\bm{f}^{(2,2)}+\cdots, \\
& \cdots
\end{split}\tag{Eq (Hierarchical)}
\end{equation}
where, $\bm{x}^{(i)}$, $\bm{p}^{(i)}$, and $\bm{f}^{(i)}$ respectively denote the input signals, parameters or variables, and contribution to the remainder of the equation of subcircuit $i$ of the given circuit. The overlapping part between $\bm{x}^{(i)}$ and $\bm{f}^{(i)}$ depends on the common nodes of the subcircuits. As shown in \ref{eq:hierarchical-equation}, $\bm{f}^{(i)}$ may be further decomposed into $\bm{f}^{(i,1)},\bm{f}^{(i,2)},\cdots$ as required.

Note that \ref{eq:flat-equation} represents only transient (TRAN) equations, and \ref{eq:flat-equation} is usually discretized in the time direction during numerical solution. At each time step, the system of algebraic equations is solved using the Newton-Raphson method
\cite[Section 7.1]{fijnvandraat2002time}
\[
\text{ Solve }\bm{x},\text{ Subject to }
\frac{1}{\beta\Delta t}\bm{Q}(\bm{x},\bm{p})+F(\bm{x},\bm{p})+\bm{b}=\bm{0},
\]
where $\bm{Q},\bm{F}$, and sparse Jacobian matrix $\nabla_{\bm{x}}\bm{Q},\nabla_{\bm{x}}\bm{F}$ need to be repeatedly calculated. For other types of simulation such as DC analysis and AC small signal analysis, \ref{eq:flat-equation} must be converted (Appendix \ref{appendix:TRAN-to-AC-equation}). Because the processing of each analysis equation is similar, the following uses TRAN analysis and a simple JSON netlist subcircuit definition (Code \ref{lst:size-dependent-resistor}) as an example to describe how to denote the calculation of $\bm{Q}^{(\cdots)},\bm{F}^{(\cdots)}$, $\nabla_{\bm{x}^{(\cdots)}\text{ or }\bm{p}^{(\cdots)}}\bm{Q}^{(\cdots)}$, and $\nabla_{\bm{x}^{(\cdots)}\text{ or }\bm{p}^{(\cdots)}}\bm{F}^{(\cdots)}$ in hierarchical circuit simulation as the forward and backward pass of a computational graph (Figure \ref{fig:computational-graph}).

\subsection{Subcircuit Module Definition in JSON Format Netlist}
\label{subsec:subckt-module-definition}
Similar to Verilog-AMS \cite[Section 6]{verilog2014verilog}, we define a circuit module that contains five parts of information (Table \ref{tab:subckt-module-definition}): (1) external nodes; (2) internal nodes; (3) input parameters; (4) decomposition of internal subcircuits; and (5) intrinsic parameters.
\begin{lstlisting}[language=json,basicstyle=\small,numbers=none,
caption={User-defined subcircuit named SizeDepResistor in the netlist: a resistor whose resistance is size-dependent},
label=lst:size-dependent-resistor]
"SizeDepResistor":{ # Define the subcircuit module.
  "ExternalNodes":["l","r"],
  "InputParams":["Rlength","Rwidth"],
  "InternalNodes":[],
  "SubModel":{
    "Expr":"[1e2*Rlength/Rwidth,]",
    "IntrinsicParams":["RValue"]
  },
  "Schematic":{
    # Instantiate each subcircuit or element in the module.
    "instanceR":{
      "MasterName":"resistor",
      "ExternalNodes":{"left":"l","right":"r"},
      "InputParams":{"resistance":"RValue"}
    }
  }
}
\end{lstlisting}
\begin{table}[htbp]
\centering
\caption{Subcircuit module definition}\label{tab:subckt-module-definition}
\begin{tabular}{l|l|l|l}
	\hline
	& \multicolumn{1}{c|}{Content} & \multicolumn{1}{c|}{Field} &\\
	\hline
	\multirow{4}{*}{\begin{tabular}[c]{@{}r@{}}Structural information\\{\small(dictionary+list)}\end{tabular}}
	& List of external node names & ExternalNodes & Required \\
	& List of internal node names &InternalNodes & Required \\
	& List of input parameter names &InputParams & Required \\
	& Internal subcircuit decomposition & Schematic & Required \\
	\hline
	\begin{tabular}[c]{@{}r@{}}Behavioral information\\{\small(differentiable function)}\end{tabular}
      & \begin{tabular}[c]{@{}l@{}}Submodel for calculating\\intrinsic parameters\end{tabular}& SubModel      & Optional \\
	\hline
\end{tabular}
\end{table}
The "Schematic" field represents internal subcircuit decomposition, which includes zero or more instantiation statements of subcircuits/devices. Each instantiation statement in "Schematic" is composed of (1) an instance name; (2) a class name (or master name); (3) external node connections; and (4) input parameter values. Code \ref{lst:size-dependent-resistor} provides an example, where the subcircuit decomposition part involves only one instance.
\begin{itemize}[partopsep=0pt,topsep=0pt,itemsep=0pt,parsep=0pt]
\item "instanceR" indicates the instance name.
\item "MasterName" indicates that the instance is of the "resistor" class. The class, or master, can be a subcircuit module or a type of built-in supported basic element.
\item "ExternalNodes" indicates that the two external nodes "left" and "right" of the instance are respectively connected to ports "l" and "r" of the module, i.e., "SizeDepResistor" here. In general case, the nodes connected to each instance in "Schematic" must come from "ExternalNodes" and "InternalNodes" of the given module.
\item "InputParams" indicates that the parameter of the instance is the internal variable "RValue" calculated by "SubModel". The parameters referenced by the instance in "Schematic" must come from:
(1) global variables; (2) "InputParams" of the module; (3) "IntrinsicParams" under "SubModel" (if any) of the model.
\end{itemize}
For more information about SubModel and its functionality, see Section \ref{subsec:EvalCompositeSubCkt}.

\subsection{Representation of Subcircuit Module Instances in a Program}
\label{subsec:subckt-instance-data-structure}
\begin{figure}[htpb]
\centering
\includegraphics[width=0.8\textwidth]{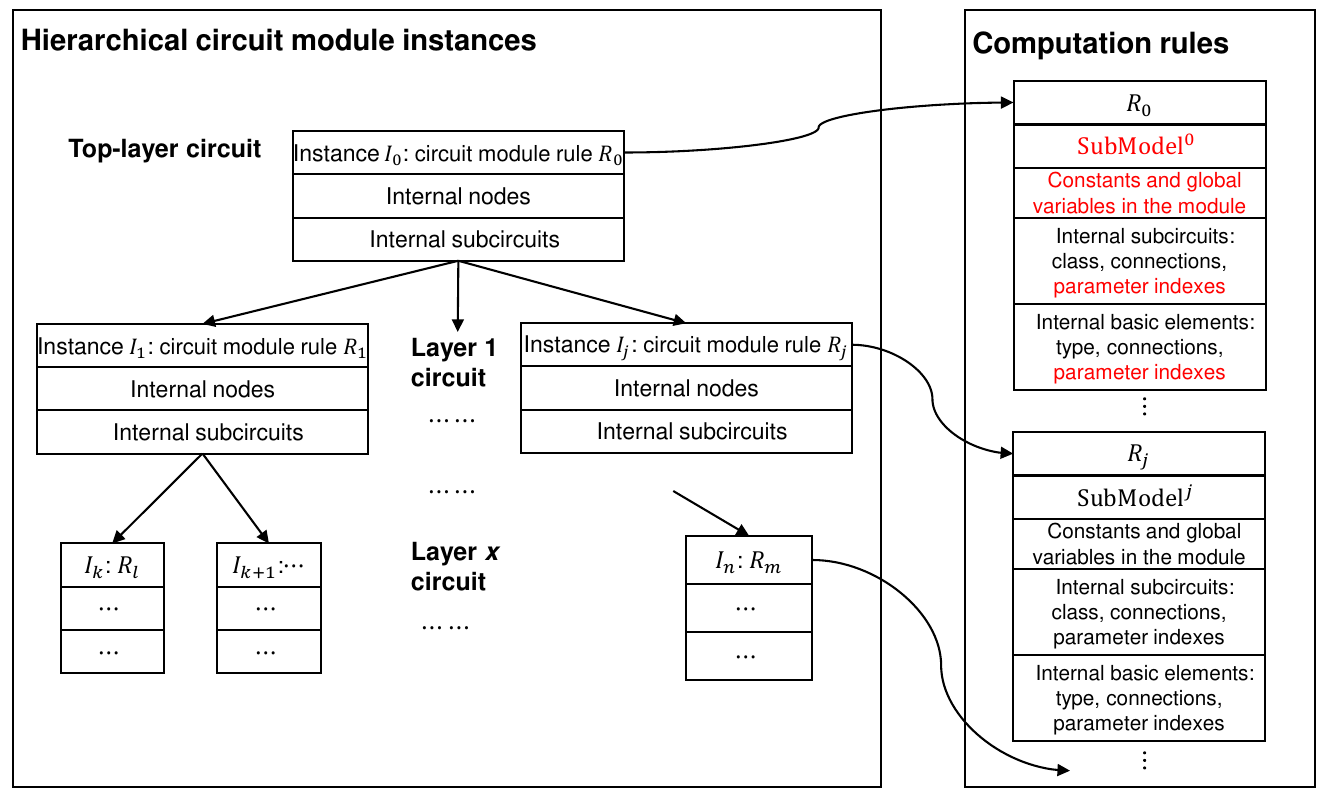}
\caption{Schematic diagram of hierarchical subcircuit module instances and computation rules. {\color{red}Text in red} indicates the parts that mark the differences from the existing method \cite{tcherniaev2003transistor}.
  Because the existing method does not need to support dynamic parameters, it is only necessary to store the fixed parameters of the devices in the computation rules of each circuit module. However, because the proposed method requires that parameters passed to lower-layer instances and devices be calculated during runtime, parameter indexing is necessary.}
\label{fig:subckt-instance-data-structure}
\end{figure}

The subcircuit definition should be compiled into an appropriate hierarchical data structure (Figure \ref{fig:subckt-instance-data-structure}) so that the equations system constructor can efficiently invoke subcircuit modules. A compiled subcircuit module contains two parts: common computation rules of subcircuits of the same master (Table \ref{tab:BasicCompositeSubCktRule}) and instance private data (Table \ref{tab:CompositeSubCkt})
\begin{table}[htbp]
\centering
\caption{Instance private data of a subcircuit module}\label{tab:CompositeSubCkt}
\begin{tabular}{l|l}
	\hline
	\multicolumn{1}{c|}{Symbol} & \multicolumn{1}{c}{Description}\\
	\hline
	rule        & Pointer to the corresponding computation rules (Table \ref{tab:BasicCompositeSubCktRule}) \\
    \textbf{in} & Internal nodes \\
	subckts     & Pointers to lower-layer subcircuit instances \\
	\hline
\end{tabular}
\end{table}
\begin{table}[htbp]
\centering
\caption{Computation rule of a subcircuit module}\label{tab:BasicCompositeSubCktRule}
\begin{tabular}{l|l}
	\hline
	\multicolumn{1}{c|}{Symbol} & \multicolumn{1}{c}{Description}\\
	\hline
	\textbf{c}       & Constants\\
	\textbf{gv}      & Global variables\\
	SubModel         & SubModel for calculating intrinsic parameters \\
	SubCktsInfo      & Lower-layer subcircuit nodes and parameter indexes \\
	BasicElementInfo & Basic element nodes and parameter indexes \\
	\hline
\end{tabular}
\end{table}
There are a few points to note:
\begin{enumerate}[partopsep=0pt,topsep=0pt,itemsep=0pt,parsep=0pt]
\item The external nodes of a subcircuit are from its upper-layer subcircuits. The top-layer circuit is a closed system without external nodes.
\item Subcircuit instances may share external nodes with one another, but the internal nodes of a subcircuit instance are exclusive to itself. When instantiating a subcircuit, ensure that its internal nodes do not conflict with each other.
\item Global variables and system signals $\bm{x}$ are globally visible to all subcircuit modules. Only the indexes \textbf{gv} of global variables need to be stored in the module's computation rules (Table \ref{tab:BasicCompositeSubCktRule}). The internal and external nodes of a module are also indexed for easy storage and passing.
\item If interactive analysis and debugging require more support information, the subcircuit master names, internal and external node parameter names, and lower-layer subcircuit instantiation statements can be added to a computation rule (Table \ref{tab:BasicCompositeSubCktRule}). Additionally, the input parameters can be dynamically recorded in an instance (Table \ref{tab:CompositeSubCkt}).
\end{enumerate}

\subsection{Instance Compilation from Subcircuit Module Definition}
\label{subsec:subckt-module-compilation}
\begin{figure}[htpb]
\centering
\begin{subfigure}{0.59\textwidth}
	\includegraphics[width = \textwidth]{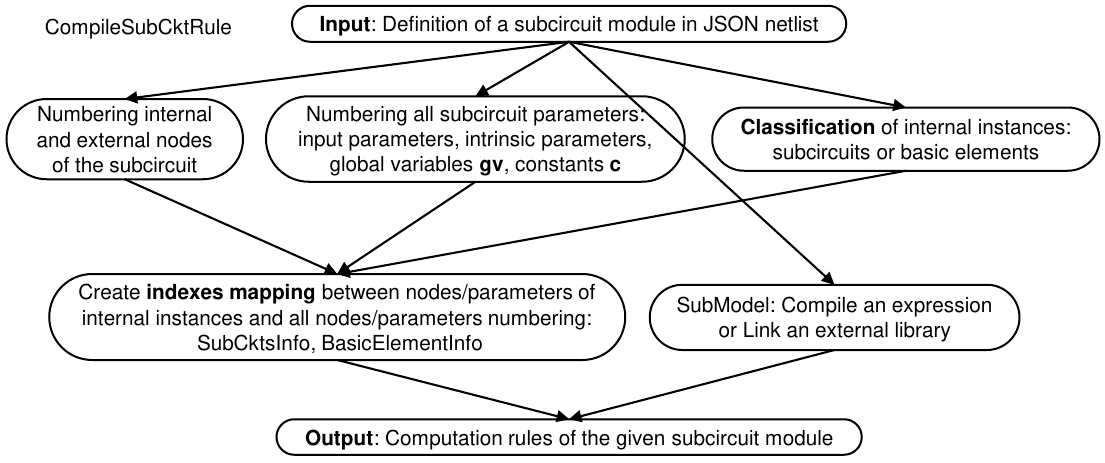}
	\caption{Compiling the computation rules of a module}
	\label{fig:compile-subckt-rule}
\end{subfigure}
\begin{subfigure}{0.35\textwidth}
	\includegraphics[width = \textwidth]{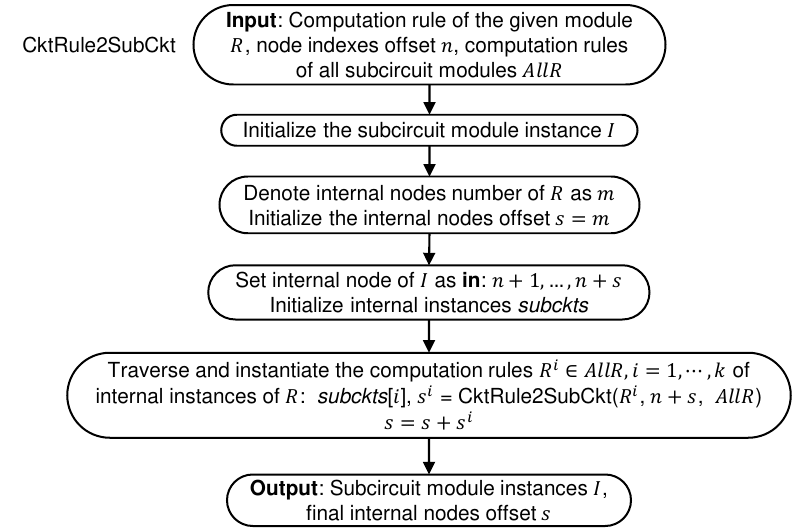}
	\caption{Recursive instantiation}
	\label{fig:cktrule-to-subckt}
\end{subfigure}
\caption{Compiling circuit modules}
\label{fig:subckt-module-compilation}
\end{figure}

A JSON netlist file can be parsed using JSON parser tools available in a variety of programming languages. The compilation of the subcircuit module definition (Section \ref{subsec:subckt-module-definition}) involves two steps:
\begin{enumerate}[partopsep=0pt,topsep=0pt,itemsep=0pt,parsep=0pt]
  \item Compile the computation rules of all subcircuit modules (Figure \ref{fig:compile-subckt-rule}). SubModel finishes parsing and compilation based on the compiler's implementation. To process other structural information (i.e., to create nodes and parameters indexes), only basic algorithms and data structures such as lists and dictionaries are needed.
\item Recursively instantiate the hierarchical circuit modules (Figure \ref{fig:cktrule-to-subckt}).
  The indexes of the input nodes are offset by $n=0$ if the instantiation program was launched at top layer circuit. The proposed method ensures that the internal nodes of each subcircuit module are independent of each other.
\end{enumerate}
Note that the internal subcircuit decomposition of a circuit module may include both basic elements and other circuit modules defined in the netlist. As such, the compiler should be able to distinguish between these two types of instances before it creates indexes for nodes and parameters. This is in addition to the compiler being able to check that there are no circular definitions of subcircuit class, undefined subcircuit modules, disconnected subcircuits, and unused nodes in the circuit.
\paragraph{Basic Elements}
\addcontentsline{toc}{subsubsection}{\ \ \ \ Basic Elements}
Basic elements are the smallest grained subcircuits without internal nodes or devices. Table \ref{tab:basic-elements-partial-list} provides a brief list of supported basic elements. To add a type of basic elements, we need to define the electrical response function for each analysis and then provide information such as external nodes and input parameters to the compiler.
\begin{table}[htpb]
\centering
\caption{Supported basic elements}
\label{tab:basic-elements-partial-list}
\begin{tabular}{l|l|l|l}
	\hline
	\multicolumn{1}{c}{MasterName}   & \multicolumn{1}{|c|}{ExternalNodes} &
	\multicolumn{1}{c|}{InputParams} & Remarks             \\
	\hline
	resistor   & left,right              & resistance  & Resistor           \\
	capacitor  & input,output            & capacitance & Capacitor           \\
	inductor   & input,output            & inductance  & Inductor           \\
	CS         & input,output            & current     & Current Source           \\
	VS         & input,output            & voltage     & Voltage Source           \\
	VCCS       & left,right,input,output & MF          & Voltage Controlled Current Source   \\
	CCCS       & iorigin,input,output    & MF          & Current Controlled Current Source   \\
	VCVS       & left,right,input,output & MF          & Voltage Controlled Voltage Source   \\
	CCVS       & iorigin,input,output    & MF          & Current Controlled Voltage Source   \\
	\hline
\end{tabular}
\end{table}

According to the modified nodal analysis method \cite{ho1975modified}, the basic elements of the voltage source type must take the branch current as one of the degrees of freedom --- this branch current is processed as an external GALV node in the compiler. At compile time, the GALV node needs to be added to the upper-layer module as an internal node. A generalized external GALV node can also be added for basic elements of a non-voltage source type such as resistors, to indicate the current flowing through the element branch. Consider a TRAN analysis example, with the resistance of the resistor denoted as $R$, the left and right nodes $l$ and $r$, the voltage value $x_l$ and $x_r$, and with or without the GALV node and current value $i,x_i$, we can present the remainder of the equation corresponding to the resistor using a sparse vector as follows:
\begin{description}
\item[Without external GALV:] $\bm{Q}=\bm{0},\bm{F}=[(l,-\frac{x_l-x_r}{R}),(r,\frac{x_l-x_r}{R})]$
\item[With external GALV:] $\bm{Q}=\bm{0},\bm{F}=[(l,-x_i),(r,x_i),(i,x_r-x_l+R\cdot x_i)]$
\end{description}
The two equations correspond to the so-called "element stamps" of the same type of elements described in \cite[Section 2.4.4]{najm2010circuit}. We may also consider them as two network analysis methods \cite{ho1975modified,hachtel1971sparse} for the same type of elements, which requires the support of both the compiler and the equations system constructor. For different analysis type, the calculation of remainder terms and that of the gradients in the basic element simulation equation must be distinguished --- we will not discuss that in detail here.

\subsection{Execution: Forward and Backward Pass of a Computational Graph}\label{subsec:EvalCompositeSubCkt}
Each basic computing unit of the computational graph (Figure \ref{fig:computational-graph}) corresponds to a subcircuit instance. When the subcircuit is invoked in a computational graph, the computing unit first takes external nodes and input variables as input from upper layer circuit. The compute unit then traverses the internal subcircuit and basic devices to calculate the equation remainder and the signal and variable gradients. Finally, these results are returned to the upper layer. See Algorithm \ref{alg:EvalCompositeSubCkt} for the internal process details.

\begin{figure}[htpb]
\centering
\includegraphics[width=0.7\textwidth]{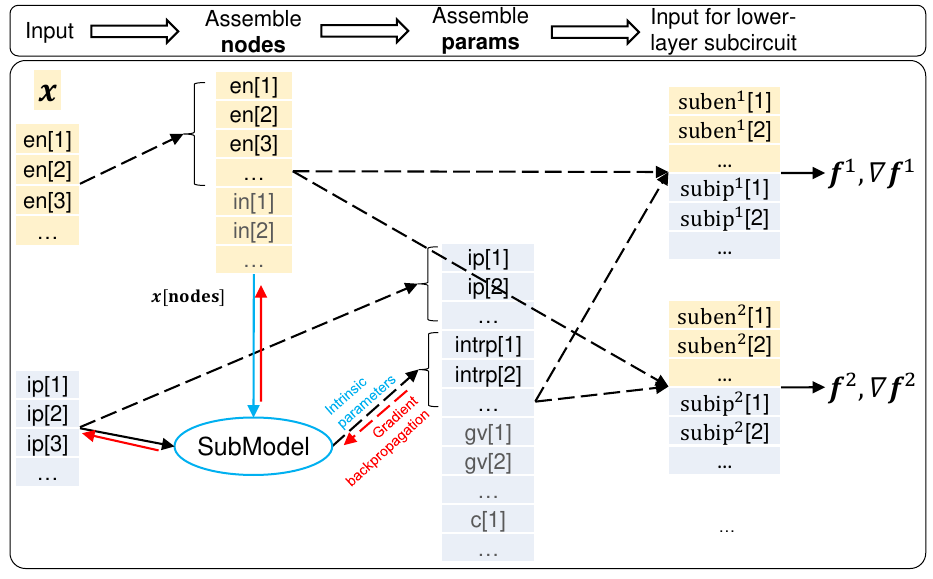}
\caption{Schematic diagram of steps 1 to 5 of Algorithm \ref{alg:EvalCompositeSubCkt}, also a zoom-in view of the calls to the top-layer to layer-1 subcircuits in Figure \ref{fig:computational-graph}. In the existing method, the input and intrinsic parameters are computed at compile time, and no gradient of parameters are propagated backward at runtime.}
\label{fig:EvalCompositeSubCkt}
\end{figure}

Figure \ref{fig:EvalCompositeSubCkt} shows steps 1 to 5 of Algorithm \ref{alg:EvalCompositeSubCkt}, where \textbf{en}, \textbf{ip}, \textbf{in}, \textbf{gv}, and \textbf{intrp} stand for external nodes, input parameters, internal nodes, global variables, and intrinsic parameters, respectively. The internal and external nodes \textbf{en},\textbf{in} of the circuit may be used to index the generalized signal values $\bm{x}[\textbf{en}],\bm{x}[\textbf{in}]$, respectively. The variables/parameters $\bm{p}$ involved in the circuit module consists of four parts: \textbf{ip}, \textbf{gv}, \textbf{intrp}, and \textbf{c}.

\begin{algorithm}[h]
\caption{Calling a Subcircuit\\
    equations remainder, signal gradient, variable gradient = \\
    {\color{white}\tiny PLACEHOLDER}EvalCompositeSubCkt($\bm{x}$,ckt,\textbf{en},\textbf{ip})
    }
\label{alg:EvalCompositeSubCkt}
\SetAlgoLined
Input: System signals $\bm{x}$, subcircuit instance ckt,
  external node indexes \textbf{en}, input parameters \textbf{ip}\;
\# Internal information of ckt: Internal nodes \textbf{in}, SubModel, global variables \textbf{gv}, constants \textbf{c}\;
1. Assemble the internal nodes \textbf{in} and external nodes \textbf{en} of ckt, resulting in \textbf{nodes}=[\textbf{en},\textbf{in}]\;
2. Calculate the intrinsic parameters according to the internal and external signals and input parameters: \textbf{intrp} = SubModel($\bm{x}$[\textbf{nodes}],\textbf{ip})\;
3. Assemble all variables and parameters of ckt, resulting in \textbf{params}= [\textbf{ip},\textbf{intrp},\textbf{gv},\textbf{c}]\;
4. Extract the external nodes \textbf{suben}$\subset$\textbf{nodes} and input parameters \textbf{subip}$\subset$\textbf{params} of each subcircuit (subckt) in ckt from \textbf{nodes},\textbf{params}, and call EvalCompositeSubCkt($\bm{x}$,subckt,\textbf{suben},\textbf{subip})\;
5. Extract the external nodes and input parameters of each basic element in ckt from \textbf{nodes},\textbf{params}, and calculate the equation remainder and gradient of each basic element\;
6. Collect the remainder terms of all equations in steps 4 and 5\;
7. Propagate signal and variable gradients of lower-layer subcircuits and basic elements backward according to the index mapping of steps 1 to 5\;
Output: Equation remainder, signal gradient, variable gradient
\end{algorithm}

\paragraph{SubModel and Intrinsic Parameters}
\addcontentsline{toc}{subsubsection}{\ \ \ \ SubModel and Internal Parameters}
In the computational graph, SubModel takes the internal and external signals and input parameters of the module as inputs, and outputs all \textit{intrinsic parameters} via \textbf{intrp}=SubModel(\textbf{signals},\textbf{ip}), where \textbf{nodes}=[\textbf{en},\textbf{in}], \textbf{signals}=$\bm{x}[\textbf{nodes}]$,
which can be passed to the lower-layer subcircuits and basic elements.
This setup is based on the Assumption \ref{assumption:intrinsic-params-dependencies}: "The behavior of a circuit module is and only is determined by the internal and external signals as well as input variables." Thus, for a SubModel, there is no need to perceive the internal signals of lower-layer subcircuits nor the node signals or parameters of other irrelevant modules. This setting can cover a considerably wide range of nonlinear effects and is easy to program.
\begin{assumption}\label{assumption:intrinsic-params-dependencies}
The intrinsic parameters in a subcircuit module are uniquely determined by the bias signals of the internal and external nodes and the input parameters of the module.
\end{assumption}
The submodel in the circuit definition should provide sufficient information so that the compiler can register the SubModel with the common rules (Table \ref{tab:BasicCompositeSubCktRule}).
In addition, there should be some protocol between the computational graph and the SubModel for obtaining the Jacobian matrix of \textbf{intrp} with respect to \textbf{signals},\textbf{ip}
\begin{equation}\label{eq:submodel-jacobian}
J_{\textbf{s}} = \nabla_{\textbf{signals}}\textbf{intrp},
J_{\textbf{ip}}=\nabla_{\textbf{ip}}\textbf{intrp},
\end{equation}
The specific implementation depends on the programming language used. As such, details are not provided here.

\paragraph{Layer-wise Gradient Backpropagation}
\addcontentsline{toc}{subsubsection}{\ \ \ \ Layer-wise Gradient Backpropagation}
A computational graph completes the computation process by calling subcircuits. The computing logic differs from that involved in the existing method (Figure \ref{fig:equations-system-constructor}) in one major aspect: In the computational graph, the input parameters \textbf{subip} of lower-layer modules or elements come from a subset of the assembled parameters \textbf{params} (Figure \ref{fig:EvalCompositeSubCkt}). Consequently, gradient backpropagation for \textbf{subip} is required. The following describes the gradient backpropagation process in Algorithm \ref{alg:EvalCompositeSubCkt} for TRAN simulation as an example.

For TRAN analysis, the returned value of Algorithm \ref{alg:EvalCompositeSubCkt} contains the following eight items:
$\bm{Q}, \bm{F}, \nabla_{\bm{x}}\bm{Q}, \nabla_{\bm{x}}\bm{F}$,
$\nabla_{\textbf{gv}}\bm{Q}, \nabla_{\textbf{ip}}\bm{Q}$,
$\nabla_{\textbf{gv}}\bm{F}$, and $\nabla_{\textbf{ip}}\bm{F}$.
Because the gradient backpropagation of $\bm{Q}$ is the same as that of $\bm{F}$, only $\bm{Q}$ is considered for simplicity.
The computation results of all subcircuits of Algorithm \ref{alg:EvalCompositeSubCkt} are recorded as $\{\bm{Q}^{i}\}$, $\{\nabla_{\bm{x}}\bm{Q}^{i}\}$, $\{\nabla_{\textbf{gv}}\bm{Q}^{i}\}$, and $\{\nabla_{\textbf{subip}^i}\bm{Q}^{i}\}$, where the superscript $i$ denotes the sequence number of the internal subcircuit or basic element. $\bm{Q}^{i}$, $\nabla_{\bm{x}}\bm{Q}^{i}$, and $\nabla_{\textbf{gv}}\bm{Q}^{i}$ can be directly assembled as
\[
\bm{Q} = \sum_i \bm{Q}^{i},
\nabla_{\bm{x}}\bm{Q} = \sum_i \nabla_{\bm{x}}\bm{Q}^{i},
\nabla_{\textbf{gv}}\bm{Q} = \sum_i \nabla_{\textbf{gv}}\bm{Q}^{i},
\]
while the gradient backpropagation of $\nabla_{\textbf{subip}^i}\bm{Q}^{i}$ needs to be processed differently based on the index of $\textbf{subip}^i$ to \textbf{params}=[\textbf{ip},\textbf{intrp},\textbf{gv},\textbf{c}].

\begin{enumerate}[partopsep=0pt,itemsep=0pt,parsep=0pt]
\item If $\text{subip}^i[j]\in\textbf{c}$, backpropagation is not performed.
\item If $\text{subip}^i[j]\in\textbf{ip}\cup\textbf{gv}$, then $\nabla_{\text{subip}^i[j]}\bm{Q}^i$ is propagated to the corresponding $\nabla_\textbf{ip}\bm{Q}$ or $\nabla_\textbf{gv}\bm{Q}$.
\item If $\text{subip}^i[j]=\textbf{intrp}[l]$ for any index $l$ (Figure \ref{fig:EvalCompositeSubCkt}), given Assumption \ref{assumption:intrinsic-params-dependencies} and the Jacobian matrix of the intrinsic parameters with respect to the signal and input parameters (Equation \ref{eq:submodel-jacobian}), then (let $\bm{g}\triangleq\nabla_{\textbf{subip}^i[j]}\bm{Q}^{i}$)
\begin{equation}\label{eq:intrinsic-params-backward}
\nabla_{x[\textbf{nodes}]}\bm{Q}\mathrel{+}=J_{\textbf{s}}[:,l]\otimes\bm{g},
\nabla_{\textbf{ip}}\bm{Q}\mathrel{+}=J_{\textbf{ip}}[:,l]\otimes\bm{g}.
\end{equation}
where, $\otimes$ represents the outer product of two vectors.
\end{enumerate}

\section{Applications}\label{sec:applications}
\subsection{CMOS Device Model: Equivalent Circuit Decomposition + Dynamic Parameters}\label{subsec:cmos-model}
\begin{figure}[htpb]
  \centering
  \includegraphics[width=0.7\textwidth]{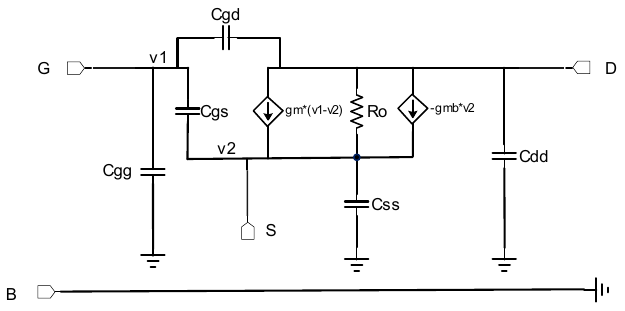}
  \caption{CMOS equivalent small-signal model \cite[Figure 2.39]{razavi2002design}, where, Ro is a resistor with resistance $\frac{1}{\text{GDS}}$.}
  \label{fig:mos-small-signal-model-2}
\end{figure}

As mentioned earlier, any subcircuit module (e.g., CMOS) that satisfies Assumption \ref{assumption:intrinsic-params-dependencies} can be modeled as a submodel-based representation featuring "equivalent circuit decomposition + dynamic parameters". This section provides an implementation example based on a lookup table. Reimplementing BSIM model \cite{chauhan2012bsim} as a SubModel is a conventional method that offers better compatibility with the existing method, but it is not adopted in this work. The CMOS module definition consists of the following elements (details about the definition are provided in Appendix \hyperref[appendix:mos-subckt]{B} and the equivalent circuit diagram under AC analysis is given in Figure \ref{fig:mos-small-signal-model-2}):
\begin{enumerate}[partopsep=0pt,topsep=0pt,itemsep=0pt,parsep=0pt]
  \item Internal and external nodes: \textbf{nodes}=[gate,source,drain,bulk];
  \item Input parameters for the device size: \textbf{ip}=[MosL,MosW];
  \item Intrinsic parameters output by SubModel:
    \textbf{intrp}=[ID,GDS,CDD,CSS,CGG,CGS,CGD, GM,GMB], whose value is determined by the four bias voltage values (\textbf{nodes}) and device size (\textbf{ip}).
\end{enumerate}
The compiler loads external libraries and generates a function object of class "lut.MosLookup" to register it as a submodel in the subcircuit rule (Table \ref{tab:BasicCompositeSubCktRule}). Among the intrinsic parameters, ID indicates the DC current between source and drain under DC analysis, and GDS, CDD, GM, etc. are equivalent small-signal parameters under AC analysis.
There are a few points to note:
\begin{enumerate}[partopsep=0pt,topsep=0pt,itemsep=0pt,parsep=0pt]
  \item The role of built-in basic elements ICS and ACVCCS (Appendix \hyperref[appendix:mos-subckt]{B}) is to ensure that ID only functions under DC analysis, while GDS, GM and GMB only function under AC analysis.
  \item The DCAC hybrid analysis or DC analysis computational graph of the equation constructor can execute this circuit module. However, the pure AC analysis computational graph cannot independently run this circuit module: In order to establish the AC analysis equations, it is necessary to first compute [GDS,GM], etc., which are determined by the DC bias voltage. This is different from directly inducing small signal linear equations through TRAN analysis equations (Appendix \hyperref[appendix:TRAN-to-AC-equation]{A}). In fact, Assumption \ref{assumption:intrinsic-params-dependencies} also stipulates that internal variables can depend on the bias voltage signal, but not on the small signals in linear analysis.
  \item The SubModel can freely call external programs, such as using three-dimensional spline interpolation, provided that it ensures compliance with the interface requirements of the corresponding automatic differentiation system.
\end{enumerate}
This submodel-based device model representation method features "equivalent circuit decomposition + dynamic parameters" and shows the following advantages:

\begin{enumerate}[partopsep=0pt,topsep=0pt,itemsep=0pt,parsep=0pt]
  \item Decoupled from circuit network analysis or simulation, the submodel is only responsible for calculating the intrinsic parameters and Jacobian matrix (Section \ref{subsec:subckt-instance-data-structure}). Circuit connectivity is not the submodel's concern.
  \item The syntax and capability boundary of a submodel in calculating intrinsic parameters depend on the compiler's processing of the "SubModel" field in the netlist. This can be implemented easily using various external programs and automatic differential tools.
\end{enumerate}

\subsection{OpAmp Device Sizing: DC Operating Points Optimization Under Different PVT Combinations}
\begin{figure}[htpb]
  \centering
  \includegraphics[width=\textwidth]{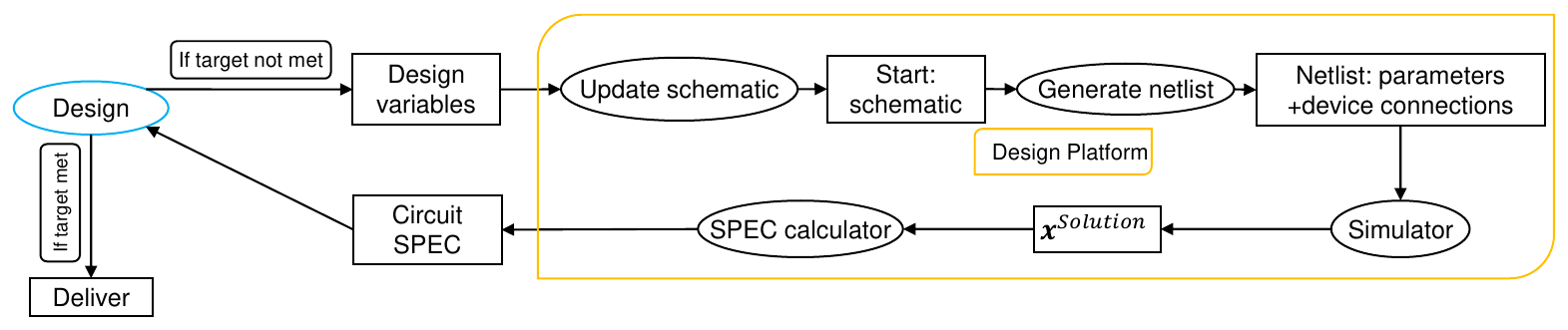}
  \caption{Manual device sizing: an iterative process}
  \label{fig:manually-design}
\end{figure}

Figure \ref{fig:manually-design} shows the device sizing process in analog circuit design. Specifically, designers connect available devices accessible to the target process into circuits, and adjust the size of each device (such as a CMOS device) based on specific methodologies and experience so that a circuit can fulfill specification requirements under a given area and power consumption constraints. In this process, repeated circuit simulation is done to quantitatively inspect the behavior and performance of a circuit without the need to manufacture the physical circuit.

\begin{figure}[htpb]
  \centering
    \includegraphics[width = 0.7\textwidth]{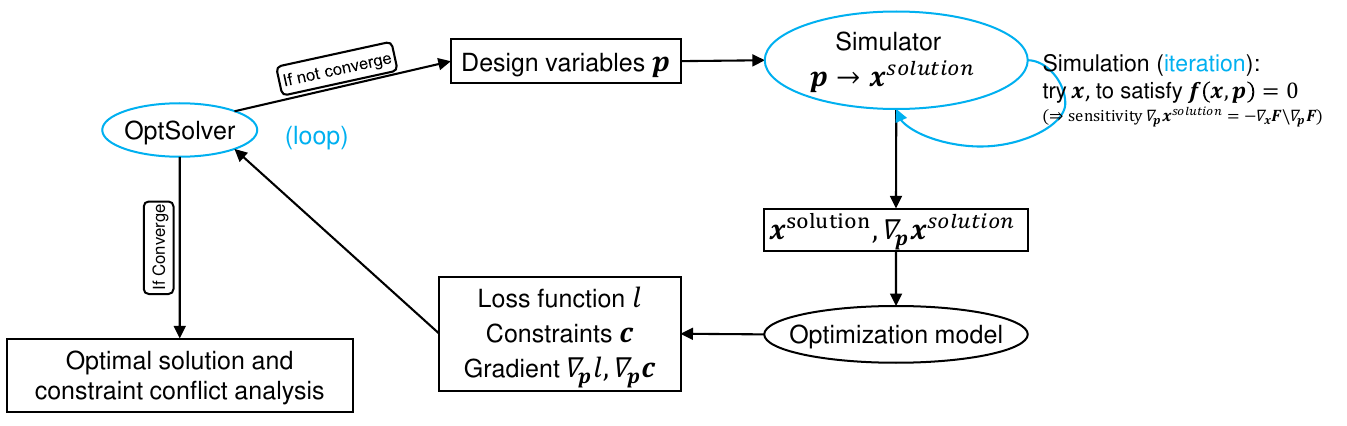}
  \caption{Automatic device sizing}
  \label{fig:solve-then-optimize}
\end{figure}

This process can be naturally converted into an optimization problem. The optimization strategy varies depending on whether the gradient is acquirable \cite{zhan2004optimization,agrawal2006circuit,huang2013efficient,nieuwoudt2005multi,peng2016efficient,girardi2011analog,lyu2018batch,wang2014enabling,lyu2017efficient,tang2018parametric}. Take DC simulation as an example. To obtain the gradient, the DC steady-state equations $\bm{F}(\bm{x},\bm{p})=\bm{0}$ naturally provide an implicit mapping from the parameter $\bm{p}$ to the solution $\bm{x}^{solution}$ of the systems of equations, with the Jacobian matrix of this mapping being $\nabla_{\bm{p}}\bm{x}^{solution}=-\nabla_{\bm{x}}\bm{F}\backslash\nabla_{\bm{p}}\bm{F}$, where $\nabla_{\bm{x}}\bm{F},\nabla_{\bm{p}}\bm{F}$ may be directly given by the equations system construction method (computational graph \ref{fig:computational-graph}) described in Section \ref{sec:Joanna}.
With this information, the gradient optimization method (Figure \ref{fig:solve-then-optimize}) can be used.
Note that in an optimization process, the inverse of $\nabla_{\bm{x}}\bm{F}$ does not need to be completely solved.
Instead, it is sufficient to solve a set of linear equations only once during each iteration's gradient backpropagation for a given loss function or constraint function $l$:
$\nabla_{\bm{p}}l(\bm{x}^{solution})=(\nabla_{\bm{p}}\bm{x})^T\cdot\nabla_{\bm{x}}l
=-(\nabla_{\bm{p}}\bm{F})^T\cdot\big(\nabla_{\bm{x}}\bm{F}^T\backslash\nabla_{\bm{x}}l\big)$

\begin{figure}[htbp]
  \centering
  \begin{subfigure}{0.3\textwidth}
    \includegraphics[width=\textwidth]{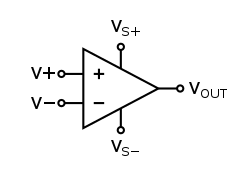}
    \caption{OpAmp schematic diagram \cite{OpAmpPNG}}
    \label{subfig:amplifier-diagram}
  \end{subfigure}
  \begin{subfigure}{0.65\textwidth}
    \includegraphics[width=\textwidth]{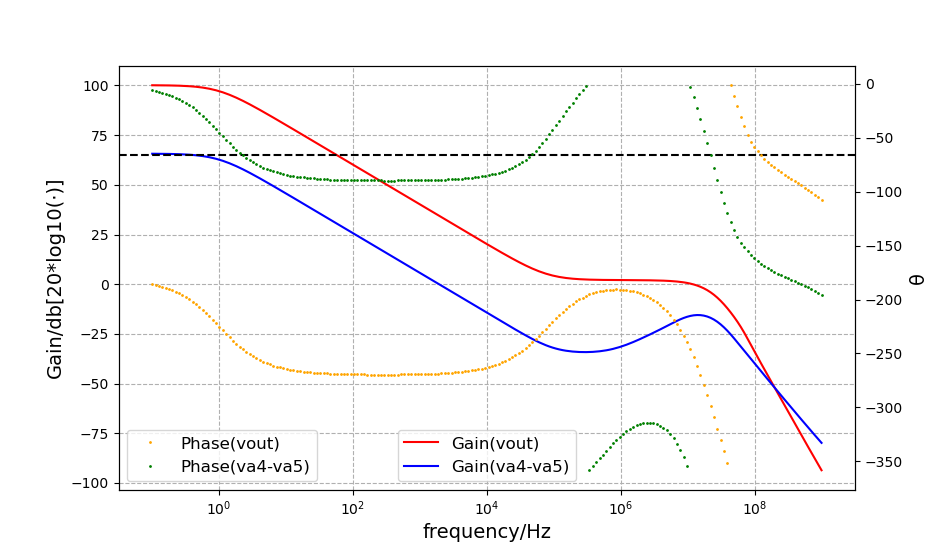}
    \caption{Input and output frequency response curves of the OpAmp after sizing}
    \label{subfig:amplifier-transfer-function}
  \end{subfigure}
  \caption{\textbf{(a)}: OpAmp schematic diagram, including the bias circuit and the main circuit, and a total of 17 n-MOSFET and 17 p-MOSFET devices. At the DC operating point, when small-signal disturbance with a given frequency is applied at $V_+,V_-$, an output signal is detected at $V_{out}$. \textbf{(b)}: Frequency response curves of OpAmp after sizing under $Corner=tt,Temperature=27$ conditions, where, va4 and va5 are the two internal nodes in the circuit.}
\end{figure}
Taking an OpAmp (Figure \ref{subfig:amplifier-diagram}) as an example, consider the design variables of the circuit (i.e., the channel length and width of each MOSFET) as the variables to be optimized.
Given the external current and voltage sources Ibias0, Ibias1, and $V_{dd},V_+,V_-$ and the load resistance and capacitance $R_L=200\Omega,C_L=10^{-10}\text{F}$, set the optimization goals as follows:
\begin{enumerate}[partopsep=0pt,topsep=0pt,itemsep=0pt,parsep=0pt]
  \item The DC operating points of all MOSFETs must be saturated under nine PVT conditions defined by $Corner\in[tt,ff,ss],Temperature\in[27,-40,125]$. For example, for an NMOS, its DC bias voltage must satisfy
    \[
      \min(V_{gs},V_{ds},V_{sb},V_{gs}-V_{th})\geq0,
    \]
    where, $V_{th}$ is subject to MosL,$V_{gs},V_{ds}$. For different PVT conditions, the SubModel needs to load different databases, and therefore the final simulation solutions obtained are also different.
  \item Under the typical condition of $Corner=tt,Temperature=27$, slight fluctuation is allowed for voltage sources $V_+,V_-$ as long as $V_++V_-=5v$ is satisfied. For the DC bias of $V_{out}$, the maximum must be greater than 4.35 V, and the minimum must be less than 0.3 V.
    Because our model (Section \ref{subsec:cmos-model}) is not used in the TRAN analysis, this requirement plays a similar role as the output swing indicator of circuits.
  \item AC analysis is performed on the circuit under the typical condition $Corner=tt$, $Temperature=27$ and a $v_{in+},v_{in-}=\pm0.5$ signal is applied at $V_+,V_-$. The DC gain $gain=20\cdot\log_{10}(|v_{out}|)$ of $V_{out}$ must reach 100.
  \item The design variables of each device meet given symmetry constraints. For example, the input MOSFET pair $M_{n0\_in},M_{n0\_ip}$ have the same size (MosL,MosW), and the current mirrors $M_{p30\_mirr},M_{p20\_mirr},M_{p10\_mirr},M_{p50\_mirr},M_{p60\_mirr}$ have the same channel length (MosL).
\end{enumerate}
\begin{equation}\label{eq:optimization}
  \begin{split}
    \min_{\bm{p}} l &= \max(5-\log_{10}(|\bm{v}[out]|),0)^2 \\
    \st \ \ \ \ \ \ 
    & \forall c\in[tt,ff,ss],t\in[27,-40,125], \\
    & \ \ \ \ \bm{x}_L\preceq\bm{x}^{c,t}\preceq\bm{x}_U;
    Saturation(\bm{x}^{c,t},\bm{p})\succeq\bm{0}; \\
    & \bm{x}^{down}[out]\leq0.3;\bm{x}^{up}[out]\geq4.35; \\
    & \bm{v}=\bm{A}^{tt,27}\backslash\bm{b}^{tt,27}; \ C\cdot\bm{p}=\bm{0}.
  \end{split}
\end{equation}
This design task can be expressed as a constrained optimization problem Prob.\eqref{eq:optimization}. $\bm{p}\to\{\bm{x}^{c,t}\},\bm{x}^{down},\bm{x}^{up}$ is obtained by solving the system of DC equations under corresponding PVT conditions with input bias $V_+,V_-$. And $\bm{v}=\bm{A}^{tt,27}\backslash\bm{b}^{tt,27}\triangleq\bm{A}\ \backslash\bm{b}$ solves the system of AC linear equations under the $Corner=tt$,$temperature=27$ condition, where the matrix elements of $\bm{A}$ are GM, GDS, etc. of each device subject to $\bm{x}^{tt,27},\bm{p}$. We can use a computational graph of mixed DCAC analysis to calculate $\bm{A},\bm{b},\nabla_{\bm{x}}\bm{A},\nabla_{\bm{x}}\bm{b}$ \footnote{$A=i\omega\cdot\nabla_{\bm{x}}\bm{Q}+\nabla_{\bm{x}}\bm{F}$. For a simpler graph implementation, use the DCAC computational graph (instead of $\nabla\bm{Q},\nabla\bm{F}$) to calculate $\nabla A$. This avoids calculating and propagating backward the second derivative of $\bm{Q},\bm{F}$ with regard to $\bm{x}$.}, which can further enable the calculation of $\nabla_{\bm{x}}l,\nabla_{\bm{p}}l$ (Appendix \hyperref[appendix:inv-linear-equation-grad]{C}). $C\cdot\bm{p}=\bm{0}$ represents the direct constraints on design variables, such as a symmetry constraint.

The optimization algorithm is implemented by using the open-source software Ipopt \cite{wachter2006implementation} and includes 72 variables, 27 equality constraints, and 308 inequality constraints to be solved. It took 356 seconds to run the whole process (including compiling Julia code and parsing netlist, etc.) on six threads on Intel(R) Core(TM) i7-8700 CPU at 3.20 GHz. Figure \ref{subfig:amplifier-transfer-function} shows the frequency response curve of the optimized circuit. The experimental results show that:
\begin{enumerate}[partopsep=0pt,topsep=0pt,itemsep=0pt,parsep=0pt]
  \item In hierarchical circuit simulation or sizing based on the computational graph, the device model and solution algorithm are decoupled from each other, allowing for high flexibility and efficiency.
  \item The parameters in the computational graph are processed to function as variables, making gradient optimization of many indicators simpler and easier.
\end{enumerate}
Note that the preceding experiments only consider the operation points of devices and circuit DC gains under typical conditions. To complete the design, more indicators (even discrete value indicators) need to be introduced into the optimization problem. There is no shortcut to properly integrating all indicators into the optimization framework, which however will not be addressed here.

\section{Conclusion}
In this paper, the static parameters of the circuit are processed as runtime variables in simulation, and the structural information and behavioral information of the circuit module/device are decoupled as "equivalent subcircuit decomposition + submodel-computed dynamic parameters". These further derive the computational graph representation of the equations system constructor for hierarchical circuits with circuit modules as the compute units of the computational graph. According to the two simple examples, this approach facilitates the decoupling and flexible interaction between netlists, models, and simulation and optimization algorithms. However, some problems exist with this approach. For example, because the variable gradient will be passed across the layers of a circuit, the topology analysis for circuit equations solvability and DAE-Index no longer works, requiring a more general hierarchical analysis theory. In the future, this approach will gain better generalization by supporting BSIM and more simulation types with more effects (i.e., S parameter or thermal effect) considered. Faster simulation is also possible if the program itself is optimized and the support for fast-SPICE technology is added.

\section*{Acknowledgment}
\addcontentsline{toc}{section}{Acknowledgment}
I would like to thank my colleagues at HiSilicon. It's been both a pleasure and an honor to work with them. This work has been also supported by many EDA and HDL experts, including Waisum Wang, Jiewen Fan, Zuochang Ye, Yang Lu, Shangxia Fang, Guoyong Shi, and Yu Ji. Special thanks go to Zhenya Zhou, Jiawei Zhuang, Long Wang, and Yuwei Fan for their kind help.

\addcontentsline{toc}{section}{References}
\bibliographystyle{unsrtnat}
\bibliography{ref}
\section*{Appendix}
\addcontentsline{toc}{section}{Appendix}
\renewcommand\thesubsection{\Alph{subsection}}
\subsection{TRAN Analysis Induced AC Analysis Equation}\label{appendix:TRAN-to-AC-equation}
To solve the transient equation \ref{eq:flat-equation} and the DC steady-state equation $\bm{F}(\bm{x}, \bm{p})=\bm{0}$,
the Newton-Raphson method calculates $\bm{Q},\bm{F}\in\mr^N$ and the Jacobian matrix $\nabla_{\bm{x}}\bm{Q},\nabla_{\bm{x}}\bm{F}\in\mr^{N\times N}$ repeatedly.
In AC small-signal analysis, consider applying perturbation $\delta\bm{x},\delta\bm{p}$ around the steady-state solution $\bm{x},\bm{p}$ and linearizing \ref{eq:flat-equation} as follows:

\[
\nabla_{\bm{x}}\bm{Q}\cdot\dot{\delta}\bm{x}
+\nabla_{\bm{p}}\bm{Q}\cdot\dot{\delta}\bm{p}
+\nabla_{\bm{x}}\bm{F}\cdot\delta\bm{x}
+\nabla_{\bm{p}}\bm{F}\cdot\delta\bm{p}=\bm{0}
\]
Let $\delta\bm{x},\delta\bm{p}$ be a small signal with an angular frequency of $\omega$: 
$\delta\bm{x}=\bm{\epsilon}_{\bm{x}}\cdot e^{i\omega t},\delta\bm{p}=\bm{\epsilon}_{\bm{p}}\cdot e^{i\omega t}$.
Then the system of linear equations of AC small-signal analysis is
\begin{equation}\label{eq:flat-ac-equation}
(i\omega\cdot\nabla_{\bm{x}}\bm{Q}+\nabla_{\bm{x}}\bm{F})\cdot\bm{\epsilon}_{\bm{x}}
= (i\omega\cdot\nabla_{\bm{p}}\bm{Q}+\nabla_{\bm{p}}\bm{F})\cdot\bm{\epsilon}_{\bm{p}}
\end{equation}

\subsection{CMOS Subcircuit}\label{appendix:mos-subckt}
\begin{lstlisting}[language=json,numbers=none,
caption={CMOS subcircuit},label=lst:mos-subckt]
"NMOSTYPE":{
  "ExternalNodes":["gate","source","drain","bulk"],
  "InputParams":["MosL","MosW"],
  "InternalNodes":[],
  "SubModel":{
    "Analysis":["DC","TRAN"],
    "ModelLoader":"SimInfo->lut.MosLookup(\"NMOSTYPE\",
      /path/to/data; SimInfo=SimInfo)",
    "IntrinsicParams":
      ["ID","GDS","CDD","CSS","CGG","CGS","CGD","GM","GMB"]
  },
  "Schematic":{
    "ids":{
      "MasterName":"ICS",
      "ExternalNodes":{"input":"source","output":"drain"},
      "InputParams":{"dc":"ID","ac":0}
    },
    "template":{
      "MasterName":"MosSmallSignalTemplate",
      "ExternalNodes":{
        "gate":"gate","source":"source",
        "drain":"drain","bulk":"bulk"
      },
      "InputParams":{
        "GDS":"GDS","CDD":"CDD","CSS":"CSS","CGG":"CGG",
        "CGS":"CGS","CGD":"CGD","GM":"GM","GMB":"GMB"
      }
    }
  }
}
\end{lstlisting}
\begin{lstlisting}[language=json,numbers=none,
caption={MosSmallSignalTemplate: Small signal equivalent circuit decomposition},label=lst:mos-small-signal-subckt]
"MosSmallSignalTemplate":{
  "ExternalNodes":["gate","source","drain","bulk"],
  "InputParams":["GDS","CDD","CSS","CGG","CGS","CGD","GM","GMB"],
  "InternalNodes":[],
  "Schematic":{
    "infr":{
      "MasterName":"resistor",
      "ExternalNodes":{"left":"drain","right":"source"},
      "InputParams":{"resistance":1e1000}
    },
    "gds":{
      "MasterName":"ACVCCS",
      "ExternalNodes":{"left":"drain","right":"source","input":"drain","output":"source"},
      "InputParams":{"MF":"GDS"}
    },
    "cdd":{
      "MasterName":"capacitor",
      "ExternalNodes":{"input":"drain","output":"bulk"},
      "InputParams":{"capacitance":"CDD"}
    },
    "css":{
      "MasterName":"capacitor",
      "ExternalNodes":{"input":"source","output":"bulk"},
      "InputParams":{"capacitance":"CSS"}
    },
    "cgg":{
      "MasterName":"capacitor",
      "ExternalNodes":{"input":"gate","output":"bulk"},
      "InputParams":{"capacitance":"CGG"}
    },
    "cgs":{
      "MasterName":"capacitor",
      "ExternalNodes":{"input":"gate","output":"source"},
      "InputParams":{"capacitance":"CGS"}
    },
    "cgd":{
      "MasterName":"capacitor",
      "ExternalNodes":{"input":"gate","output":"drain"},
      "InputParams":{"capacitance":"CGD"}
    },
    "gm":{
      "MasterName":"ACVCCS",
      "ExternalNodes":{
        "left":"gate","right":"source","input":"drain","output":"source"
      },
      "InputParams":{"MF":"GM"}
    },
    "gmb":{
      "MasterName":"ACVCCS",
      "ExternalNodes":{
        "left":"bulk","right":"source","input":"drain","output":"source"
      },
      "InputParams":{"MF":"GMB"}
    }
  }
}
\end{lstlisting}

\subsection{Gradient Backpropagation of the Solution of a Linear Equations System}\label{appendix:inv-linear-equation-grad}
Consider a system of real linear equations $\bm{A}(\bm{x})\bm{v}=\bm{b}(\bm{x})$ with respect to $\bm{v}$,
where $\bm{A},\bm{b}$ are nonlinearly dependent on $\bm{x}$ being a sparse matrix/vector, and
$\nabla_{\bm{x}}\bm{A},\nabla_{\bm{x}}\bm{b}$
are computable, then $\bm{v}$ will also be nonlinearly dependent on $\bm{x}$. Differentiating the system of equations yields:
\[
(\bm{A}+\nabla_{\bm{x}}\bm{A}\cdot\ud\bm{x})\cdot
(\bm{v}+\nabla_{\bm{x}}\bm{v}\cdot\ud\bm{x})
=\bm{b}+\nabla_{\bm{x}}\bm{b}\cdot\ud\bm{x}
\]
When the zeroth and second order terms are dropped, the following equation holds for any $\ud\bm{x}$
\begin{equation}\label{eq:inv-full-differentiation}
\nabla_{\bm{x}}\bm{A}\cdot\ud\bm{x}\cdot\bm{v}
+\bm{A}\cdot\nabla_{\bm{x}}\bm{v}\cdot\ud\bm{x}
=\nabla_{\bm{x}}\bm{b}\cdot\ud\bm{x},
\end{equation}
Assume a loss function $l(\bm{v})$, whose gradient $\nabla_{\bm{v}}l$ can be calculated.
Backpropagating the gradient to $\bm{x}$ is equivalent to finding the solution of
$\nabla_{\bm{x}}l=\nabla_{\bm{v}}l\cdot\nabla_{\bm{x}}\bm{v}$.
In fact, we do not need to actually calculate $\nabla_{\bm{x}}\bm{v}$ and store the data, which can be rather dense.
According to Equation \ref{eq:inv-full-differentiation}, we have
\[
\begin{split}
& \nabla_{\bm{x}}\bm{v}\cdot\ud\bm{x}
= \bm{A}^{-1}\cdot\big(\nabla_{\bm{x}}\bm{b}\cdot\ud\bm{x}-
\nabla_{\bm{x}}\bm{A}\cdot\ud\bm{x}\cdot\bm{v}\big),\forall \ud\bm{x}, \\
\Rightarrow
& \nabla_{\bm{v}}l\cdot\nabla_{\bm{x}}\bm{v}\cdot\ud\bm{x}
= (\nabla_{\bm{v}}l\cdot\bm{A}^{-1})\cdot\big(\nabla_{\bm{x}}\bm{b}\cdot\ud\bm{x}-
\nabla_{\bm{x}}\bm{A}\cdot\ud\bm{x}\cdot\bm{v}\big),\forall \ud\bm{x}, \\
\end{split}
\]
Therefore, in order to calculate $\nabla_{\bm{v}}l\cdot\nabla_{\bm{x}}\bm{v}$,
we only need to solve the sparse matrix linear equations system $\nabla_{\bm{v}}l\cdot\bm{A}^{-1}$ once in advance,
element-wisely set the values of $\ud\bm{x}$ to 1 and the rest to 0. Then, we obtain
$\nabla_{\bm{x}}l=\nabla_{\bm{v}}l\cdot\nabla_{\bm{x}}\bm{v}$.

For the case of a set of complex linear equations, a similar discussion can also be carried out using the Wirtinger derivative.

\end{document}